\documentclass[a4paper,12pt]{amsart}
\usepackage{amssymb}
\usepackage{ifthen}
\usepackage{cite}
 \usepackage[dvips]{graphicx}
\nonstopmode\numberwithin{equation}{section}
\usepackage{amssymb}
\usepackage{ifthen}
\usepackage{graphicx}
\usepackage{amsmath}
\usepackage[T1]{fontenc} 
\usepackage[utf8]{inputenc}
\usepackage[usenames,dvipsnames]{color}
\usepackage{color}
\usepackage[english]{babel}
\usepackage{fancyhdr}
\usepackage{fancybox}
\usepackage{tikz}

\newtheorem*{theoA}{Theorem A}
\newtheorem*{theoB}{Theorem B}
\newtheorem*{theoC}{Theorem C}
\newtheorem*{theoD}{Theorem D}
\newtheorem*{theoE}{Theorem E}
\newtheorem*{theoF}{Theorem F}

\nonstopmode \numberwithin{equation}{section}
\theoremstyle{plain}

\newtheorem{conj}{Conjecture}

\theoremstyle{definition}
\newtheorem{defi}{Definition}[section]

\newtheorem{thm}{Theorem}[section]
\newtheorem{prob}{Problem}[section]
\newtheorem{cor}{Corollary}[section]
\newtheorem{ques}{Question}[section]
\newtheorem{prop}{Proposition}[section]
\newtheorem{rem}{Remark}[section]
\newtheorem{lem}{Lemma}[section]
\newtheorem{note}{Note}[section]

\newtheorem{Open Ques}{Open question}[section]

\newcounter{minutes}
\divide\time by 60
\newcounter{hours}
\multiply\time by 60
\addtocounter{minutes}{-\time}

\newcounter {own}
\def\theown {\thesection       .\arabic{own}}

\newenvironment{pf}[1][]{%
 \vskip 3mm
 \noindent
 \ifthenelse{\equal{#1}{}}%
  {{\slshape Proof. }}%
  {{\slshape #1.} }%
 }%
{\qed\bigskip}

\newcounter{alphabet}

\def\be{\begin{equation}}
\def\ee{\end{equation}}

\newcommand{\bee}{\begin{enumerate}}
\newcommand{\eee}{\end{enumerate}}

\newcommand{\blem}{\begin{lem}}
\newcommand{\elem}{\end{lem}}
\newcommand{\bthm}{\begin{thm}}
\newcommand{\ethm}{\end{thm}}
\newcommand{\bcor}{\begin{cor}}
\newcommand{\ecor}{\end{cor}}
\newcommand{\beg}{\begin{examp}}
\newcommand{\eeg}{\end{examp}}
\newcommand{\begs}{\begin{examples}}
\newcommand{\eegs}{\end{examples}}
\usepackage{mathrsfs}
\newcommand{\bdefn}{\begin{defn}}
\newcommand{\edefn}{\end{defn}}

\newcommand{\bprob}{\begin{prob}}
\newcommand{\eprob}{\end{prob}}
\newcommand{\bei}{\begin{itemize}}
\newcommand{\eei}{\end{itemize}}

\newcommand{\bcon}{\begin{conj}}
\newcommand{\econ}{\end{conj}}
\newcommand{\bcons}{\begin{conjs}}
\newcommand{\econs}{\end{conjs}}
\newcommand{\bprop}{\begin{prop}}
\newcommand{\eprop}{\end{prop}}
\newcommand{\br}{\begin{rem}}
\newcommand{\er}{\end{rem}}
\newcommand{\brs}{\begin{rems}}
\newcommand{\ers}{\end{rems}}
\newcommand{\bo}{\begin{obser}}
\newcommand{\eo}{\end{obser}}
\newcommand{\bos}{\begin{obsers}}
\newcommand{\eos}{\end{obsers}}
\newcommand{\bpf}{\begin{pf}}
\newcommand{\epf}{\end{pf}} 
\newcommand{\ba}{\begin{array}}
\newcommand{\ea}{\end{array}}
\newcommand{\beq}{\begin{eqnarray}}
\newcommand{\beqq}{\begin{eqnarray*}}
\newcommand{\eeq}{\end{eqnarray}}
\newcommand{\eeqq}{\end{eqnarray*}}

\begin{document}

\title{Improved and Refined Bohr-Type Inequalities for Slice Regular Functions over Octonions}

\author{Sabir Ahammed}
\address{Sabir Ahammed, Department of Mathematics, Jadavpur University, Kolkata-700032, West Bengal, India.}
\email{sabira.math.rs@jadavpuruniversity.in}

\author{Molla Basir Ahamed$^*$}
\address{Molla Basir Ahamed, Department of Mathematics, Jadavpur University, Kolkata-700032, West Bengal, India.}
\email{mbahamed.math@jadavpuruniversity.in}

\subjclass[{AMS} Subject Classification:]{ Primary 30G35. Secondary 30B10, 17A35, 30H05. }

\def\thefootnote{}
\footnotetext{ {\tiny File:~\jobname.tex,
printed: \number\year-\number\month-\number\day,
          \thehours.\ifnum\theminutes<10{0}\fi\theminutes }
} \makeatletter\def\thefootnote{\@arabic\c@footnote}\makeatother
\maketitle
\pagestyle{myheadings}
\markboth{S. Ahammed and M. B. Ahamed }{Bohr phenomena for slice regular functions over octonions}
\begin{abstract}
A crucial extension of quaternionic function theory to octonions is the concept of slice regular functions, introduced to handle holomorphic-like properties in a non-associative setting. In this paper, first we present a generalization of the Bohr inequality, and improved versions of the Bohr inequality for slice regular functions over the largest alternative division algebras of octonions $\mathbb{O}$. Moreover, we provide a refined version of the Bohr inequality for slice regular functions $f$ on $\mathbb{B}$ such that $  {\rm Re}(f(x)) \leq 1 $ for all $x \in \mathbb{B}$. All the results are shown to be sharp.
\end{abstract}
\section{\bf Introduction}
In the study of Dirichlet series, Bohr's remarkable work on power series in complex analysis has generated significant research activity in complex analysis and related areas. Let $\mathcal{H}(\mathbb{D})$ be denotes the class of analytic functions $f(z)=\sum_{k=0}^{\infty}a_kz^k$  in the open unit disk $\mathbb{D}:=\{ z\in \mathbb{C}:|z|<1\}$ such that $|f(z)|\leq 1$ for all $z\in \mathbb{D}.$  H. Bohr \cite{Bohr-1914} discovered the following interesting phenomenon for class $\mathcal{H}(\mathbb{D})$.
\begin{theoA}\cite{Bohr-1914} If $f(z)=\sum_{k=0}^{\infty}a_kz^k\in \mathcal{H}(\mathbb{D}),$  then 
	\begin{align}\label{Eq-1.11}
		\sum_{k=0}^{\infty}|a_k||z|^k\leq 1 \;\;\;\mbox{for}\;\;\; |z|\leq \dfrac{1}{3}.
	\end{align}
	The constant $1/3$ is best possible.
\end{theoA}
The inequality $\eqref{Eq-1.11}$ fails when $ |z|>{1}/{3} $ in the sense that there are functions in $ \mathcal{H}\left(\mathbb{D}\right) $ (e.g., $f_a(z)=(a-z)/(1-az)$, where $a\in\mathbb{D}$) for which the inequality is reversed when $ |z|>{1}/{3} $. H. Bohr initially showed in \cite{Bohr-1914} that the inequality \eqref{Eq-1.11} holds only for $|z|\leq 1/6$, which was later improved independently by M. Riesz, I. Schur, F. Wiener and some others. The sharp constant $1/3$  and the inequality \eqref{Eq-1.11} in Theorem A are famously known as, respectively, the Bohr radius and the classical Bohr inequality for the family $ \mathcal{H}\left(\mathbb{D}\right) $.  Several other proofs of this interesting inequality were given in different articles  (see \cite{Sidon-1927,Paulsen-Popescu-Singh-PLMS-2002,Tomic-1962}). \vspace{1.2mm}

In $1995,$ the Bohr inequality was applied by Dixon \cite{Dixon-BLMS-1995} to the long-standing problem of characterizing Banach algebras that satisfying the von Neumann inequality and then this application garnered considerable attention from various researchers, leading to investigations of the phenomenon in a variety of function spaces. However, it can be noted that not every class of functions has Bohr phenomenon, for example, B\'en\'eteau \emph{et al.} \cite{Beneteau- Dahlner-Khavinson-CMFT-2004} showed that there is no Bohr phenomenon in the Hardy space $ H^p(\mathbb{D},X), $ where $p\in [1,\infty)$ and Hamada \cite{Hamada-Math.Nachr.-2021} showed that there is no Bohr-radius for holomorphic mappings with values in the unit ball of a complex Hilbert space ${H}$ with $\dim {H}\geq 2.$  Moreover, $\mathcal{H}\left(\mathbb{D}\right)$ is not the only class of analytic functions where the Bohr radii are studied; many other classes of functions and some integral operators also have their Bohr radii examined.  The study of the Bohr phenomenon for holomorphic functions in several complex variables and in Banach spaces has becomes very active (see \cite{Aizenberg-PAMS-2000, Galicer-Mansilla-Muro-TAMS-2020, Kumar-PAMS-2022, Lata-Singh-PAMS-2022, Liu-Ponnusamy-PAMS-2021, Hamada-RM-2024}); in the non-commutative setting ( see \cite{Paulsen-Singh-BLMS-2006,Popescu-IEOT-2019,Paulsen-Popescu-Singh-PLMS-2002}), for the operator-valued generalizations of the Bohr theorem in the single-variable case, (see \cite{ Popescu-TAMS-2007}) and references therein. For other aspects of the Bohr radius problem, we refer to the articles \cite{Kayumov-Khammatova-Ponnusamy-CRMA-2020,Allu-Halder-PEMS-2023,Hamada-Math.Nachr.-2021,Blasco-CM-2017,Hamada-IJM-2009,Hamada-AAMP-2025,Arora-Kumar-Ponnusamy-MM-2025,Kumar-Manna-JMAA-2023,Kumar-Manna-BSM-2025} and references cited therein.

\subsection{The slice regular functions and the Bohr inequality} The octonions, denoted by $\mathbb{O}$ form an $8$-dimensional normed division algebra over $\mathbb{R}$. They generalize quaternions but differ significantly due to the loss of associativity (octonions are only alternative). Developing a function theory over octonions presents unique challenges due to their non-commutativity and non-associativity. A crucial extension of quaternionic function theory to octonions is the concept of slice regular functions, introduced to handle holomorphic-like properties in a non-associative setting. Slice regular functions over octonions naturally extend holomorphic function theory to non-associative settings. While they share properties with quaternions, challenges arise due to the loss of associativity. Their study is closely linked to modern research in mathematical physics, geometry, and functional analysis. Octonions play a role in string theory, supergravity, and gauge theories, and slice regularity is connected to $G_2$-manifolds and special holonomy spaces. The slice regular functions are used in solving PDEs, generalizing Laplace and Dirac equations in $8$-dimension, and defining spectral theory and operator calculus. Despite their success, challenges remain, including defining differentiation in non-associative algebras, extending the Cauchy integral formula, and developing an octonionic residue theorem.\vspace{2mm}

As a generalization of holomorphic functions of one complex variable, the theory of slice regular functions of one quaternionic variable was initiated by Gentili and Struppa \cite{Gentili-Struppa-AVM-2007} and further developed by Clifford algebras \cite{Colombo-Sabadini-Struppa-IJM-2009} and octonions \cite{Gentili-Struppa-RMJM-2010}. Based on the concept of stem functions, these three classes of functions were eventually unified and generalized into real alternative algebras \cite{Ghiloni-Perotti-AdvM-2011}. Following the historical path, the Bohr inequality has recently been generalized to the non-commutative (but associative) algebra of quaternions for slice regular functions in \cite{Rocchetta-Gentili-Sarfatti-MN-2012}.  Rocchetta \emph{et al.} \cite{Rocchetta-Gentili-Sarfatti-MN-2012} have established a special case of Theorem A for quaternions, where regular quaternionic rational transformations from \cite{Bisi-Gentili-IUMJ-2009,Stoppato-AG} are heavily utilized. This is because quaternions form a skew but associative field. When restricted to each slice, the function \( f \) described in Theorem B can be viewed, according to \cite[Lemma 3.2]{Xu-PRSE-2021}, as a vector-valued holomorphic function mapping from \( \mathbb{D} \) into the unit ball of \( \mathbb{C}^4 \) endowed with the standard Euclidean norm. Blasco \cite[Theorem 1.2]{Blasco-OTAA-2010} introduced the Bohr radius for holomorphic functions mapping from $\mathbb{D}$ into the unit ball of $\mathbb{C}^n$ (with $n \geq 2$) and demonstrated that it is zero. From this perspective, the theory of slice regular functions differs from that of vector-valued holomorphic functions.\vspace{2mm}

In \cite{Xu-PRSE-2021}, Xu established the Bohr inequality and its sharpness for slice regular functions over the non-commutative and non-associative algebra of octonions, which is the largest (finite-dimensional) alternative division algebras. 
\begin{theoB}\emph{\cite[Theorem 1.3]{Xu-PRSE-2021}}
	Let $f(x)=\sum_{k=0}^{\infty}x^ka_k$ with $a_k\in \mathbb{O}$ be a slice regular function in the open unit ball $\mathbb{B}$ of $\mathbb{O}$ such that $|f(x)|\leq 1$ for all $x\in \mathbb{B}.$ Then 
	\begin{align*}
		\sum_{k=0}^{+\infty}|x^ka_k|\leq 1,\; |x|\leq \dfrac{1}{3}.
	\end{align*}
Moreover, the constant $1/3$ is sharp.
\end{theoB}

Paulsen \emph{et al.} \cite{Paulsen-Popescu-Singh-PLMS-2002} were the first to demonstrate that the Bohr radius depends on the power of the initial term $|f(0)|$ of the majorant series. They proved that the Bohr radius is $1/2$ when $|f(0)|$ is replaced by $|f(0)|^m$, where $m=2$. Later, Blasco \cite{Blasco-OTAA-2010} extended the range of $m \in [1,2]$ and established a refinement of Theorem A on $\mathbb{D}$.

\begin{theoC}\emph{\cite[Proposition 1.4]{Blasco-OTAA-2010}}
	If $ f(z)=\sum_{k=0}^{\infty}a_kz^k\in\mathcal{H}(\mathbb{D}) $
	and $m\in (1,2]$, then
	\begin{align*}
		|a_0|^m+\sum_{k=1}^{\infty}|a_k|r^k
		\leq 1 \;\; \mbox{for}\;\; |z|=r\leq \frac{m}{2+m}.
	\end{align*}
	The number $m/(2+m)$ cannot be improved.
\end{theoC} 
The study of generalizations, improvements, and sharper versions of Bohr-type inequalities has been an active research topic in recent years. Many researchers continuously investigated refined Bohr-type inequalities and also examining their sharpness for certain classes of analytic functions. For detailed information on such studies, the readers are referred to the articles  
\cite{Evdoridis-Ponnusamy-Rasila-RM-2021,Liu-Ponnusamy-PAMS-2021,Liu-Ponnusamy-Wang-RACSAM-2020,Liu-Liu-Ponnusamy-BSM-2021,Ponnusamy-Vijayakumar-Wirths-JMAA-2022,Allu-Arora-JMAA-2022} and the references therein.\vspace{2mm} 

 Inspired by the results in \cite{Blasco-OTAA-2010}, this paper generalizes Theorem B, offering an analogue of Theorem C. Furthermore, our proof technique deviates from existing approaches. 
\begin{thm}\label{Thm-1.4}
	Let $m\in (0,2]$ and $f(x)=\sum_{k=0}^{\infty}x^ka_k$ with $a_k\in \mathbb{O}$ be a slice regular function in the open unit ball $\mathbb{B}$ of $\mathbb{O}$ such that $|f(x)|\leq 1$ for all $x\in \mathbb{B}.$ Then 
	\begin{align*}
		\mathcal{A}_f(x):=|a_0|^m+\sum_{k=1}^{\infty}|x^ka_k|\leq 1,\;\;\mbox{for}\;\; |x| \leq R_{m},
	\end{align*}
	where $R_{m}:={m}/{(2+m)}.$ The constant $R_{m}$ is best possible for each $m$.
\end{thm}
\begin{rem}
	Setting $m=1$, Theorem \ref{Thm-1.4} reduces exactly to Theorem B, confirming that Theorem \ref{Thm-1.4} extends the results of Theorem B. Moreover, Theorem \ref{Thm-1.4} remains valid for $ m\in (0, 2]$ thereby providing a analogue of Theorem C, which was only applicable for $ m\in (1, 2]$.
\end{rem}
Refining Bohr inequality for certain classes of functions is an important part of the Bohr phenomenon. Kayumov and Ponnusamy \cite{Kayumov-Ponnusamy-CRACAD-2018} have studied Theorem A for its refined version and established the following sharp result for functions in the class $\mathcal{H}(\mathbb{D})$.
\begin{theoD}\emph{\cite[Theorem 3]{Kayumov-Ponnusamy-CRACAD-2018}}
	If $ f(z)=\sum_{k=0}^{\infty}a_kz^k\in\mathcal{H}(\mathbb{D})$, then 
	\begin{align*}
		\sum_{k=0}^{\infty}|a_k|r^k+|f(z)-a_0|^2 &\leq 1 \;\; \mbox{for}\;\; |z|=r\leq {1}/{3}.
	\end{align*}
	The constant $1/3$ cannot be improved.
\end{theoD}

In the next result, by using arguments similar to those in the proof of Theorem \ref{Thm-1.4}, we show that a refined version of Theorem B of the Bohr inequality can be obtained which gives an analogue of Theorem D. 
\begin{thm}\label{BS-thm-1.2}
	Let $m\in (0,2]$, $\lambda\geq 0$ and $q\geq 1$.
	Let $f$ be a slice regular function same as in Theorem B. Then
	\begin{align*}
		\mathcal{B}_f(x):=\mathcal{A}_f(x)+\lambda |f(x)-a_0|^q
		\leq 1 \;\; \mbox{for}\;\; |x|\leq R_{m,\lambda, q},
	\end{align*}
	where $R_{m, \lambda,q}$ is the unique root in $(0, 1)$ of the equation 
	\[ -\frac{m}{2}+\frac{r}{1-r}+\lambda\left(\frac{r}{1-r}\right)^q=0.\]
	In particular, 
	\begin{align*}
		\sum_{k=0}^{\infty}|x^ka_k|+ |f(x)-a_0|^2
		\leq 1\;\;\;\mbox{for}\;\;\;|x|\leq R_{1,1,2}:=2-\sqrt{3}.
	\end{align*}  
\end{thm}
An interesting application of the Bohr inequality for the class  $\mathcal{H}(\mathbb{D})$ was explored by Ponnusamy \emph{et al.} \cite{Ponnusamy-Vijayakumar-Wirths-RM-2020}, who investigated a modified version of the classical Bohr inequality, known as a refinement of the classical Bohr inequality and established the following refined Bohr inequality on $\mathbb{D}$.
\begin{theoE}\cite[Theorem 2]{Ponnusamy-Vijayakumar-Wirths-RM-2020}
	If $ f(z)=\sum_{k=0}^{\infty}a_kz^k\in\mathcal{H}(\mathbb{D}) $, then
	\begin{align*}
		\sum_{k=0}^{\infty}|a_k|r^k
		+\left(\frac{1}{1+|a_0|}+\frac{r}{1-r} \right)\sum_{k=1}^{\infty}|a_k|^2r^{2k}\leq 1 \;\; \mbox{for}\;\; |z|=r &\leq \frac{1}{2+|a_0|}.
	\end{align*}
	The numbers $1/(1+|a_0|)$ and $1/(2+|a_0|)$ cannot be improved.
\end{theoE}	

Following the trend of study on Bohr radius problems, it is natural to explore whether Theorem E can be extended to slice regular functions, and whether the corresponding radius can be shown to be sharp. In the next result, we obtain a refined version of Theorem \ref{Thm-1.4} which gives an analogue version of Theorem E.
\begin{thm}\label{Thm-1.5}
	Let $m\in (0,1]$ and $f(x)=\sum_{k=0}^{\infty}x^ka_k$ with $a_k\in \mathbb{O}$ be a slice regular function in the open unit ball $\mathbb{B}$ of $\mathbb{O}$ such that $|f(x)|\leq 1$ for all $x\in \mathbb{B}.$ Then 
	\begin{align*}
		\mathcal{C}_f(x):=\mathcal{A}_f(x)+\left(\dfrac{1}{1+|a_0|}+\dfrac{|x|}{1-|x|}\right)\sum_{k=1}^{\infty}|x^ka_k|^2\leq 1 \;\;\mbox{for}\;\; |x|\leq R_m= \dfrac{m}{2+m}.
	\end{align*}  
	The constant $R_m$ is best possible for each $m$.
\end{thm}
\subsection{The Bohr inequality in terms of the planar integral} Let $f$ be holomorphic in $\mathbb D$, and for $0<r<1$,  let $\mathbb D_r=\{z\in \mathbb C: |z|<r\}$.
Let $S_r:=S_r(f)$ denote the planar integral
\begin{align*}
	S_r=\int_{\mathbb D_r} |f^{\prime}(z)|^2 d A(z).
\end{align*}
If the function $f\in \mathcal{H}(\mathbb{D})$ has Taylor's series expansion $f(z)=\sum_{k=0}^{\infty}a_kz^k $, then we obtain  (see \cite{Kayumov-Ponnusamy-CRACAD-2018})
\begin{align*}
	S_r= \pi\sum_{k=1}^\infty k|a_k|^2 r^{2k}.
\end{align*}

In the study of the improved Bohr inequality, the quantity $ S_r $ plays a significant role. There are many results on the improved Bohr inequality for the class $ \mathcal{H}(\mathbb{D}) $ (see \cite{Ismagilov-Kayumov-Ponnusamy-JMAA-2020,Kayumov-Ponnusamy-CRACAD-2018}), and for harmonic mappings on unit disk (see \cite{Evdoridis-Ponnusamy-Rasila-IM-2019}).\vspace{2mm}

It is not merely a matter of using arbitrary combinations of \( {S_r}/{\pi} \) and its powers. The crucial point is that the inequality must hold for the given radius, which must remain fixed. Moreover, it is essential to show that this radius is optimal and that the associated coefficient bounds cannot be improved. Thus, the challenge lies in identifying the appropriate combination on the left-hand side of the inequality. The motivation for the study on Bohr phenomenon via proper combinations over octonions is based on the discussions above and the observation of the following remark concerning improved Bohr inequality.
\begin{rem}\label{Rem-2.1}
	Ismagilov \emph{et al.} \cite{Ismagilov-Kayumov-Ponnusamy-JMAA-2020} remarked that for any function $F : [0, \infty)\to [0, \infty)$ such that $F(t)>0$ for $t>0$, there exists an analytic function $f : \mathbb{D}\to\mathbb{D}$ for which the inequality
	\begin{align*}
		\sum_{n=0}^{\infty}|a_n|r^n+\frac{16}{9}\left(\frac{S_r}{\pi}\right)+\lambda\left(\frac{S_r}{\pi}\right)^2+F(S_r)\leq 1\;\; \mbox{for}\;\; r\leq\frac{1}{3}
	\end{align*}
	is false, where $\lambda$ is given in \cite[Theorem 1]{Ismagilov-Kayumov-Ponnusamy-JMAA-2020}.
\end{rem}
Considering Remark \ref{Rem-2.1}, it is noteworthy to observe that augmenting a non-negative quantity with the Bohr inequality does not yield the desired inequality for the class $\mathcal{H}\left(\mathbb{D}\right)$. This observation motivates us to study the Bohr inequality further for slice regular functions over octonions with a suitable setting.\vspace{2mm}

 It is a natural to ask the question.
\begin{ques}\label{Q-1}
	Can the improved Bohr inequality be established for slice regular functions over octonions?
\end{ques}

Before answering Question \ref{Q-1}, we need to establish some preliminaries. Accordingly, we introduce the quantity $S_{x}^*$ for slice regular functions of the form $f(x)=\sum_{k=0}^{\infty}x^ka_k$ in $\mathbb{B}$ is given by 
\begin{align}\label{Eq-1.2}
	S_{x}^*:=\sum_{k=1}^{\infty}k|x^ka_k|^2.
\end{align}
This observation motivates us to study the Bohr inequality further  for slice regular functions over octonions with a suitable setting of $S_{x}^*.$ In the next result, we first establish an improved version of the Bohr inequality  for slice regular functions over octonions.\vspace{1.2mm}

\noindent To serve our purpose, we consider a polynomial of degree $ N $ as follows
\begin{align}\label{Eq-1.3}
	Q_N(w):=d_1w+d_2w^2+\dots+d_Nw^N, \;\; \mbox{where} \;\;
	d_i\geq 0, i=1,2,\dotsm N,
\end{align}
and obtain the following improved version of the Bohr inequality.

\begin{thm}\label{BS-thm-1.3}
		Let $m\in (0,1]$ and $f(x)=\sum_{k=0}^{\infty}x^ka_k$ be a slice regular function in the open unit ball $\mathbb{B}$ of $\mathbb{O}$ such that $|f(x)|\leq 1$ for all $x\in \mathbb{B}.$ Then 
	\begin{align*}
	\mathcal{D}_f(x):=	\mathcal{A}_f(x)+Q\left(S^*_x\right)\leq 1\;\;\;\mbox{for}\;\;\;|x|\leq R_m,
	\end{align*}  
	where the coefficients of the polynomial $Q_N$ satisfy the condition
	\begin{align}\label{Eq-1.6}
		L(d_1,\dots, d_N):&=8d_1M_m^2+6c_2d_2M_m^4+\dots+2(2N-1)c_Nd_NM_m^{2N}
		\\ \nonumber
		&\leq m,
	\end{align}
	with $M_m:=m(2+m)/(4m+4)$ and 
	\[
	c_k:=\max_{x\in [0,1]}\left(x(1+x)^2(1-x^2)^{2k-2}\right),
	\quad k=2,3,\cdots, N.
	\]
	The constant $ R_m$ (for each $m$) is best possible for each $d_1,\dots, d_N$ which satisfy $\eqref{Eq-1.6}$. 
\end{thm}
\begin{note}
	Observe that the inequality \eqref{Eq-1.11} also can be expressed in terms of the distance formulation 
	\begin{align*}
		\sum_{k=1}^{\infty}|a_k||z|^k\leq 1-|a_0|=1-|f(0)|=\mbox{dist} \left(f(0), \partial\mathbb{D}\right)\;\;\mbox{for\; all }\;\; |z|\leq \dfrac{1}{3},
	\end{align*}
	where $`\mbox{dist}$' is the Euclidean distance between $f(0)$ and the boundary of $\mathbb{D}$.
\end{note}

Motivated by this formulation, researchers have extensively examined the Bohr inequality in various settings for different classes of functions (see \cite{Liu-Ponnusamy-Wang-RACSAM-2020,Ponnusamy-Vijayakumar-Wirths-JMAA-2022,Xu-Ren-JGA-2018} and references therein).  Drawing inspiration from these articles and their conclusions, Xu \cite{Xu-PRSE-2021} examined the Bohr inequality in a more generalized setting for slice regular functions over octonions and established the following result.
\begin{theoF}\cite[Theorem 1.7]{Xu-PRSE-2021}
	Let $f(x)=\sum_{k=0}^{\infty}x^ka_k$ with $a_k\in \mathbb{O}$ be a slice regular function in the open unit ball $\mathbb{B}$ of $\mathbb{O}$ such that $f(\mathbb{B})\subset \Pi:=\{ x\in \mathbb{O}: \mathrm{Re}(x)\leq 1\}.$ Then
	\[ \sum_{k=1}^{\infty}|x^ka_k|\leq \emph{dist}\left(f(0), \partial \Pi\right) \;\;\mbox{for}\;\; |x|\leq\dfrac{1}{3}.\]
\end{theoF}
In the next result, by using arguments similar to those in the proof of Theorem F and \cite[Theorem 4.2]{Xu-PRSE-2021}, we obtain a refined version of Theorem F of the Bohr inequality which gives an analogue of Theorem D. We omit the proof.
\begin{thm}\label{Th-1.5}
	Let $m\in (0,1]$, $\lambda\geq 0$ and $j\geq 1$.
	Let $f$ be a same as in Theorem F. Then
	\begin{align*}
		\mathcal{A}_f(x)+\lambda |f(x)-f(0)|^j
		\leq 1 \;\; \mbox{for}\;\; |x|\leq R^*_{m,\lambda,j},
	\end{align*}
	where $R^*_{m, \lambda,j}$ is the unique root in $(0, 1)$ of the equation 
	\[ -\frac{m}{2}+\dfrac{r}{1-r}+\lambda2^{j-1}\left(\frac{r}{1-r}\right)^j=0.\]
	In particular, we have
	 \begin{align*}
		\sum_{k=0}^{\infty}|x^ka_k|
		+|f(x)-f(0)|^2
		\leq 1 \;\;\;\mbox{for}\;\;\;|x|\leq R^*_{1,1,2}:= \sqrt{5}-2\approx 0.236068.
	\end{align*}  
\end{thm}
\begin{rem}
	For $m=1$, $\lambda=0$, we get the Theorem F from Theorem \ref{Th-1.5}. In this sense, we say that Theorem \ref{Th-1.5} is a two fold generalization of Theorem F.
\end{rem}
It is natural to ask whether an analogue of Theorem E exists when $f$ maps $\mathbb{B}$ into $\Pi.$ In this paper, we answer this question by proving the following result, subject to certain condition on $a_0=f(0)$.
\begin{thm}\label{Thm-1.7}
	Let $f$ be as in Theorem F with $ \mathrm{Re}(a_0)\in [0,1).$ Then 
	\[ \mathcal{E}_f(x):=\sum_{k=0}^{\infty}|x^ka_k|+\left(\frac{1}{1+a_0}+\frac{|x|}{1-|x|}\right)\sum_{k=1}^{\infty}|x^{k}a_k|^2\leq 1\]
	holds for all $|x|=r\leq R_*$, where $R_*\approx 0.24683$ is the unique root in $(0, 1)$ of the equation $3r^3-5r^2-3r+1=0$. The constant $R_*$ cannot be improved.
\end{thm} 
In the next theorem, we obtain an improved version of the Theorem \ref{Thm-1.7}  with a suitable setting of the quantity $S_{x}^*,$ where  $S_{x}^*$ is defined in \eqref{Eq-1.2}.
\begin{thm}\label{Theom-1.7}
 Let $f$ be as in Theorem F with $ \mathrm{Re}(a_0)\in [0,1).$ Then 
	\[ \mathcal{F}_{f}(|x|):=\sum_{k=0}^{\infty}|x^ka_k|+\left(\frac{1}{1+a_0}+\frac{|x|}{1-|x|}\right)\sum_{k=1}^{\infty}|x^ka_k|^2+\beta \left(S^*_{x}\right)\leq 1\]
	for $|x|\leq R(a_0)= 1/(5-2a_0),$ where $\beta=8/9.$ The constant $8/9$ is best possible.
\end{thm} 

\section{\bf{Preliminaries}}
In this section, we recall necessary definitions and preliminary results used in the sequel for slice regular functions from \cite{Ghiloni-Perotti-AdvM-2011}.
\subsection{The algebra of octonions}
Let $\mathbb{C},$ $\mathbb{H},$ $\mathbb{O}$ denote the algebra of complex numbers, quaternions and octonions, respectively and$\{ 1, i,j,k\}$ be the standard basis of the non-commutative, associative, real algebra of quaternions with the multiplication rules 
\begin{align*}
	i^2=j^2=k^2=ijk=-1.
\end{align*} 
For each element $a=x_0+x_1i+x_2j+x_3k$ in  $\mathbb{H}$ $(x_0,x_1,x_2,x_3\in \mathbb{R}),$ the conjugate of $a$ is defined as   $\overline{a}=x_0-x_1i-x_2j-x_3k.$ By the well-known Cayley-Dickson process, the real algebra of octonions can be built from $\mathbb{H}$ as $\mathbb{O}=\mathbb{H}+l\mathbb{H}$ with 
\begin{align*}
	&\overline{a+lb}=\overline{a}-lb, \;\; (a+lb)+(c+ld)=(a+c)+l(b+d),\;\;\\& (a+lb)(c+ld)=(ac-d\overline{b})+l(\overline{a}d+cb)\;\;\mbox{for}\;\;a,b,c,d\in \mathbb{H}.
\end{align*} 
As a consequence, $\{ 1,i,j,k,li,lj,lk\}$ form the canonical real vector basis of $\mathbb{O}$. Every element $x\in \mathbb{O}$ can be composed into the real part $\mathrm{Re}(x)=(x+\overline{x})/2$ and the imaginary part $\mathrm{Im}(x)=x-\mathrm{Re}(x).$ Define the modulus of $x$ $|x|=\sqrt{x\overline{x}},$ which is exactly the usual Euclidean norm in $\mathbb{R}^8.$ Also, the modulus is multiplicative, i.e., $|xy|=|x||y|$ for all $x,y\in \mathbb{O}.$ Every non-zero element $x\in \mathbb{O}$ has a multiplicative inverse given by $x^{-1}=|x|^{-2}\overline{x}.$ The construction above shows that $\mathbb{O}$ is a non-commutative, non-associative, normed and division algebra. See for instance \cite{Baez-BAMS-2002} for more explanation on the octonions.\vspace{1.2mm} 

The set of square roots of $-1$ in $\mathbb{O}$ is the six-dimensional unit sphere given by 
\[ \mathbb{S}:=\{q\in \mathbb{H}: q^2=-1 \}.\] 

For each $I\in \mathbb{S},$ denote by $\mathbb{C}_I:=<1,I>\cong \mathbb{C}$ the subalgebra $\mathbb{O}$ of  generated over $\mathbb{R}$ by $1$ and $I$.\vspace{1.2mm}

Notice that each $x\in \mathbb{O}$ can be expressed as $x = \alpha + \beta I_x$ with $\alpha \in \mathbb{R},$ $\beta\in \mathbb{R}^+$ and $I_x \in \mathbb{S}$. This inconspicuous observation allows decomposing $\mathbb{O}$ into `complex slices'
\[ \mathbb{O}=\bigcup_{I\in \mathbb{S}}\mathbb{C}_I,\]
which derives the remarkable notion of slice regularity over octonions.

\subsection{Slice functions}
Given an open set ${D}$ of $\mathbb{C}$, invariant under the complex conjugation, its circularization $\Omega_D$ is defined by
\[ \Omega_D=\bigcup_{I\in \mathbb{S}}\{ \alpha+\beta I: \exists\; \alpha, \beta \in \mathbb{R},\; s.t.\; z= \alpha+i \beta\in D\}.\]
A subset $\Omega$ in $\mathbb{O}$ is called to be circular if $\Omega = \Omega_D$ for some $D\subset \mathbb{C}$. The open
unit ball $ \mathbb{B}=\{ x\in \mathbb{O}:|x|<1\}$ and the right half-space $\{x\in \mathbb{O}: \mathrm{Re}(x)>0 \}$ are two
typical examples of the circular domain. For a series expansion for slice regular functions and related properties, we refer to the article \cite{Stoppato-AM-2012}.

\begin{defi}
	A function $F:D\to \mathbb{O}\otimes_\mathbb{R}\mathbb{C}$ on an open set $D\subset \mathbb{C}$ invariant
	under the complex conjugation is called a stem function if the $\mathbb{O}$-valued components $F_1,$ $F_2$ of $F=F_1+iF_2$ satisfies
	\[ F_1(\overline{z})=F_1(z),\;\; F_2(\overline{z})=-F_2(z),\;\forall\; z=\alpha+i\beta.\]
\end{defi}
Each stem function $F$ induces a (left) slice function $f = \mathcal{I}(F):\Omega_D\to \mathbb{O}$ given by
\[ f(x):=F_1(z)+IF_2(z),\;\forall\; x=\alpha+I\beta\in \Omega_D.\]
We will denote the set of all such induced slice functions on $\Omega_D$ by
\[ \mathcal{S}(\Omega_D):=\{ f=\mathcal{I}(F): F
\; \mbox{is}\; \mbox{a}\; \mbox{stem}\; \mbox{function}\; \mbox{on}\; D\}.\]
Each slice function $f$ is induced by a unique stem function $F$ since $F_1$ and $F_2$ are determined by $f$. In fact, it holds that \[ F_1(z)=\dfrac{1}{2}\left(f(x)+f(\overline{x})\right),\;\; z\in \Omega_D,\]
and \[ f(z)=\begin{cases}
	\dfrac{1}{2I_x}\left( f(x)-f(\overline{x})\right)\;\;\mbox{if}\;\;z\in \Omega_D\setminus \mathbb{R},\\\vspace{4mm}
	0,\;\;\;\;\;\;\;\;\;\;\;\;\;\;\;\;\;\;\;\;\;\;\;\;\;\;\; \mbox{if}\;\; z\in \Omega_D\cap \mathbb{R}.
\end{cases}\]
Recall that a $\mathrm{C}^1$ function $F:D\to \mathbb{O}\otimes_\mathbb{R}\mathbb{C}$  is holomorphic if, and only if, its components $F_1$, $F_2$ satisfy the Cauchy-Riemann equations 
\[ \dfrac{\partial F_1}{\partial\alpha}=\dfrac{\partial F_2}{\partial\beta},\;\;\dfrac{\partial F_1}{\partial\beta}=-\dfrac{\partial F_2}{\partial\alpha},\;\; z=\alpha+i\beta\in D. \]
\begin{defi}
	A (left) slice function $f=\mathcal{I}(F)$ on $\Omega_D$ is regular if its stem function $F$ is holomorphic on $D.$ Denote the class of slice regular functions on $\Omega_D$ by
	\[ \mathcal{SR}\left(\Omega_D\right):=\{f=\mathcal{I}(F)\in \mathcal{S}(\Omega_D): F\;\mbox{is}\;\mbox{holomorphic}\;\mbox{on}\; D\}.\]
\end{defi}
For $f\in \mathcal{SR}\left(\Omega_D\right),$ the slice derivative is defined to be the slice regular function $f^\prime$
on $\Omega_D$ obtained as
\[ f^\prime(x):= \mathcal{I}\left( \dfrac{\partial F}{\partial z}(z)\right)=\dfrac{1}{2}\mathcal{I}\left(\dfrac{\partial F}{\partial \alpha}(z)-i\dfrac{\partial F}{\partial \beta}(z) \right). \]
Recall that $\mathbb{O}$ is non-associative but alternative, \emph{i.e.}, the associator $(x,y,z):= (xy)z-x(yz)$ of three elements $x,y,z\in \mathbb{O}$ is an alternating function in its arguments. Meanwhile, the Artin theorem asserts that the subalgebra generated by two elements of $\mathbb{O}$ is associative. Hence, a class of examples of slice regular functions is given by polynomials of one octonionic variable with coefficients in $\mathbb{O}$ on the right side. Indeed, each slice regular function $f$ defined in $\mathbb{B}$ admits the expansion of convergent power series 
\[ f(x)=\sum_{k=0}^{\infty} x^ka_k,\;\;\{ a_k\}\subset \mathbb{O},\] for all $x\in \mathbb{B}.$\vspace{1.2mm}

For simplicity, let $\mathbb{B}_I$ the intersection $\mathbb{B}\cap\mathbb{C}_I$ for any $I\in \mathbb{S}$. Then the restriction $f|_{\mathbb{B}_I}$  is holomorphic on $\mathbb{B}_I$. Furthermore, the relation between slice regularity and complex holomorphy can be presented as follows.
\begin{lem}(Splitting lemma)
	Let $\{ I_0=1,I,I_1,II_1,I_2, II_2,I_3,II_3\}$ be a splitting
	basis for $\mathbb{O}.$ For $f\in \mathcal{SR}(\Omega_D),$ there exist holomorphic functions $f_n:\Omega_D\cap\mathbb{C}_I\to \mathbb{C}_I,$ $m\in \{0,1,2,3\},$ such that 
	\[ f(z)=\sum_{n=0}^{3}f_n(z)I_n,\;\; \forall\; z\in \Omega_D\cap\mathbb{C}_I.\] 
\end{lem}
Due to that the pointwise product of two slice functions is not a slice function generally, the notion of slice product was introduced.
\begin{defi}
	Let $f=\mathcal{I}(F)$ and $g=\mathcal{I}(G)$ be in $\mathcal{S}(\Omega_D)$ with stem functions $F=F_1+iF_2$ and $G=G_1+iG_2.$ Then $FG=F_1G_1-F_2G_2+i\left(F_1G_2+F_2G_1\right)$ is
	still a stem function. The slice product of $f$ and $g$ is the slice function on $\Omega_D$ defined by $f{.}g:= \mathcal{I}(FG).$
\end{defi}
In general, $f.g\neq g.f.$ If the components $F_1$, $F_2$ of the first stem function $F$ are real-valued, then $f = \mathcal{I}(F)$ is termed as slice preserving. For the slice preserving function $f$ and slice function $g$, the slice product $f.g$ coincides with $fg$.
\begin{defi}
	For $f=\mathcal{I}(F)\in \mathcal{S}(\Omega_D)$ with $F=F_1+iF_2,$ define the slice conjugate of $f$ as
	\[ f^c= \mathcal{I}\left(\overline{F_1}+i\overline{F_2}\right),\] and the normal function (or symmetrization) of $f$ as 
	\[ N(f)=f.f^c=f^c.f,\] which is slice preserving on $\Omega_D$.
\end{defi}
Let $\mathcal{Z}_f$ denote the zero set of $f$ on $\Omega_D.$ 
\begin{defi}
	Let $f\in \mathcal{S}(\Omega_D).$ If $f$ does not vanish identically, then its slice reciprocal is defined as 
	\[f^{-\bullet}(x):=N(f)(x)^{-1}.f^c(x)=N(f)(x)^{-1}f^(x) \] which is a slice function on $\Omega_D\setminus {\mathcal{Z}_{N(f)}}.$
\end{defi}
More recently, Ghiloni \emph{et. al.} \cite{Ghiloni-Perotti-AdvM-2011} found a new and nice relation between the values of reciprocals $f^{-\bullet}(x)$ and $f^{-1}(x)$ for slice functions $f\in \mathcal{S}(\Omega_D)$ as \[f^{-\bullet}(x)=f\left( T_f(x)\right)^{-1}, \] where $T_f$ is a bijective self-map of $\Omega_D\setminus E$ $E=\{ a+\beta I:I\in \mathcal{S}, z=\alpha+\beta i\in D\;\mbox{for}\; F_2(z)=0\}$ given by 
\[ T_f(x)=\left( f^c(x)^{-1}((xf^c(x)) F_2(z))\right)F_2(z)^{-1},\] which reduces to the known result $T_f(x)= f^c(x)^{-1}xf^c(x)$ for the associative
algebra of quaternions.

 \section{\bf Proof of main results}
\begin{proof}[\bf Proof of Theorem \ref{Thm-1.4}]
	Let	$f(x)=\sum_{k=0}^{\infty}x^ka_k$ with $a_k\in \mathbb{O}$ be a slice regular function in the open unit ball $\mathbb{B}$ of $\mathbb{O}$ such that $|f(x)|\leq 1$ for all $x\in \mathbb{B}.$ Under the condition of Theorem \ref{Thm-1.4}, it follows that (see \cite[Lemma 3.2]{Xu-PRSE-2021})
	\begin{align}\label{Eq-2.4}
		|a_k|\leq 1-|a_0|^2,\;\;\mbox{for}\;\;k\in \mathbb{N}.
	\end{align}
	In view of \eqref{Eq-2.4}, we obtain
	\begin{align*}
	\mathcal{A}_f(x)&\leq|a_0|^m+(1-|a_0|^2)\dfrac{|x|}{\left(1-|x|\right)}\leq 1
	\end{align*}
	if, and only if, $|x|\leq R_{m},$ where
	\begin{align*}
		R_{m}= \inf_{|a_0|<1} \{\mathcal{M}\left(|a_0|\right)\}\;\;\mbox{and}\;\;   \mathcal{M}\left(|a_0|\right):=\dfrac{1-|a_0|^m}{2-|a_0|^2-|a_0|^m}.
	\end{align*}
	Moreover, for $0<m\leq 2,$ an easy computation gives that 
	\begin{align*}
		\mathcal{M}^{\prime}(t)=\dfrac{t^{m-1}C(t)}{\left(2-t^2-t^m\right)^2},\;\;\mbox{where}\;\;C(t):=-m-(2-m)t^2+2t^{2-m}.
	\end{align*}
	Since $C^\prime(t)=2(2-m)t^{1-m}(1-t^m)\geq 0$ for $m\in (0,2]$ and $0\leq t\leq 2,$ this implies that $C(t)\leq C(1)=0$ and, hence $\mathcal{M}$ is a decreasing function of $t\in [0,1).$ Thus, it follows that
	\begin{align*}
		\mathcal{M}(t)\geq \lim\limits_{t\rightarrow 1^{-}}\mathcal{M}(t)=\dfrac{m}{2+m}=R_{m}.
	\end{align*}
	This proves the desired inequality.\vspace{1.2mm}
	
	Finally, to show the sharpness, given $a\in (0,1)$ and $u\in \partial \mathbb{B},$ we consider the slice regular function
	\begin{align}\label{Eq-1.1}
		f_a(x)=(1-xa)^{-\bullet }\ast(a-x)u=a-(1-a^2)u\sum_{k=1}^{\infty}x^ka^{k-1},\;\;x\in \mathbb{B}.
	\end{align}
	For $f_a,$ we obtain
	\begin{align*}
		\mathcal{A}_{f_a}(x)&=a^m+(1-a^2)\sum_{k=1}^{\infty}a^{k-1}|x|^k\\&=a^m+\dfrac{(1-a^2)|x|}{1-a|x|}\\&=1+(1-a)\mathcal{G}(a),
	\end{align*}
	where 
	\begin{align*}
		\mathcal{G}(a):=\dfrac{(1+a)|x|}{1-a|x|}-\left(\dfrac{1-a^m}{1-a}\right).
	\end{align*}
	Taking $a$ very closed to $1$, we see that 
	\begin{align*}
		\lim\limits_{a\rightarrow 1^{-}}\mathcal{G}(a)=\dfrac{2 |x|}{1-|x|}-m>0\; \mbox{for}\; |x|>R_m.
	\end{align*}
	Consequently, we see that $\mathcal{A}_{f_a}(x)>1$, which establishes that the constant  $R_m$ is best possible, concluding the proof. 
\end{proof}

\begin{proof}[\bf Proof of Theorem \ref{BS-thm-1.2}]
 In view of \eqref{Eq-2.4}, we obtain
 \begin{align*}
 		\mathcal{B}_f(x)&\leq |a_0|^m+(1-|a_0|^2)\sum_{k=1}^{\infty}|x|^k+\lambda \left( (1-|a_0|^2)\sum_{k=1}^{\infty}|x|^k\right)^q:=\mathcal{B}^*_f(|x|). 
 \end{align*}
For the next step in the proof, we make use of the following inequality, as established in \cite{Kayumov-Khammatova-Ponnusamy-MJM-2020}
\begin{align}\label{Eq-2.55}
	\frac{1-t^m}{1-t}\geq \dfrac{m}{2}\; \mbox{for}\; m\in (0, 2]\; \mbox{and}\; t\in [0, 1).
\end{align}
	Let $|x|=r$.  Taking into account \( |a_0| < 1 \) and the equation \eqref{Eq-2.55}, we derive that
	\begin{align*}
		\mathcal{B}^*_f(|x|)&\leq  1-(1-|a_0|^2)\dfrac{m}{2}+(1-|a_0|^2)\dfrac{r}{1-r}+\lambda (1-|a_0|^2)^q\left(\dfrac{r}{1-r}\right)^q\\&=1+(1-|a_0|^2)\left( -\dfrac{m}{2}+\dfrac{r}{1-r}+\lambda (1-|a_0|^2)^{q-1} \left(\dfrac{r}{1-r}\right)^q\right)\\&\leq 1+(1-|a_0|^2)\left( -\dfrac{m}{2}+\dfrac{r}{1-r}+\lambda \left(\dfrac{r}{1-r}\right)^q\right).
	\end{align*}
	Thus, our desired inequality $\mathcal{B}_f(x)\leq 1$ is established for $r\leq R_{m, \lambda,q},$ where $R_{m, \lambda,q}$ is the unique root in $(0, 1)$ of the equation 
	\[ -\frac{m}{2}+\frac{r}{1-r}+\lambda\left(\frac{r}{1-r}\right)^q=0.\] 
	This completes the proof.
\end{proof}

 \begin{proof}[\bf Proof of Theorem \ref{Thm-1.5}]
  In view of \eqref{Eq-2.4}, we obtain
 	\begin{align*}
 		\mathcal{C}_f(x)&\leq |a_0|^m+(1-|a_0|^2)\sum_{k=1}^{\infty}|x|^k+\left(\dfrac{1}{1+|a_0|}+\dfrac{|x|}{1-|x|}\right)(1-|a_0|^2)^2\sum_{k=1}^{\infty}|x|^{2k}\\&=1+(|a_0|^m-1)+(1-|a_0|^2)\dfrac{|x|}{1-|x|}+\left(\dfrac{1}{1+|a_0|}+\dfrac{|x|}{1-|x|}\right)\dfrac{(1-|a_0|^2)^2|x|^2}{1-|x|^2}\\&:=\mathcal{K}^*_f(|a_0|,|x|).
 	\end{align*}
 	To proceed further in the proof, we use the following inequality (see \cite{Lin-Liu-Ponnu-ACS-2021})
 	\begin{align}\label{Eq-2.5}
 		 \frac{1-t^m}{1-t}\geq m\; \mbox{for}\; m\in (0, 1]\; \mbox{and}\; t\in [0, 1).
 	\end{align}
 	Let $|a_0|=t$. In view of \eqref{Eq-2.5}, we obtain that 
 	\begin{align*}
 		\mathcal{K}^*_f(t,|x|)\leq 1+m(t-1)+K_1(1-t^2)+K_2(1-t)(1-t^2)+K_3(1-t^2)^2:=\Psi(t),
 	\end{align*}
 	where $t\in [0,1],$ and
 	\begin{align*}
 		K_1:=\dfrac{|x|}{1-|x|},\;\;K_2:=\dfrac{|x|^2}{1-|x|^2}\;\;\mbox{and}\;\; K_3:=\dfrac{|x|^3}{(1-|x|)(1-|x|^2)}.
 	\end{align*}
 	It is clear that $K_1,$ $K_2$ and $K_3$ are non-negative for all $x\in \mathbb{B}\setminus\{0\}$.\vspace{2mm} 
 	
 \noindent	Moreover, we see that 
 	\begin{align*}
 		\begin{cases}
 			\Psi^\prime(t)=m-2K_1t+K_2(3t^2-2t-1)+4K_3(t^3-t),
 			\vspace{2mm}\\
 			\Psi^{\prime\prime}(t)=-2K_1+2K_2(3t-1)+4K_3(3t^2-1),
 			\vspace{2mm}\\
 			\Psi^{\prime\prime\prime}(t)=6K_2+24tK_3.
 		\end{cases}
 	\end{align*}
 	Since $K_2$ and $K_3$ are non-negative, we have $\Psi^{\prime\prime\prime}(t)>0$ for all $t\in [0,1]$  which leads that the function $	\Psi^{\prime\prime}$ is increasing in $t,$ which implies that 
 	\begin{align*}
 		\Psi^{\prime\prime}(t)\leq \Psi^{\prime\prime}(1)&=-2K_1+4K_2+8K_3\\&= \dfrac{2|x|\left( 4|x|^2+2|x|(1-|x|)-(1-|x|^2) \right)}{(1-|x|)(1-|x|^2)}\\&=\dfrac{2|x|(1+|x|)(3|x|-1)}{(1-|x|)(1-|x|^2)}.
 	\end{align*}
 Furthermore, it is clear that \( \Psi^{\prime\prime}(t) \leq 0 \) for \( |x| \leq \frac{1}{3} \) and \( t \in [0, 1] \), which implies that \( \Psi^{\prime} \) is a decreasing function on \( [0, 1] \). Hence, a simple calculation shows that 
 	\begin{align*}
 		\Psi^{\prime}(t)\geq \Psi^{\prime}(1)=m-2K_1=\dfrac{m-(2+m)|x|}{1-|x|}\geq 0\;\;\mbox{for}\;\; |x|\leq \dfrac{m}{2+m}.
 	\end{align*}
 	Since \( \Psi^{\prime}(t) \geq 0 \) for \( t \in [0, 1] \) and \( |x| \in \left( 0, \frac{m}{m+2} \right) \), we conclude that \( \Psi \) is an increasing function on \( [0, 1] \). Thus, we obtain 
 	\begin{align*}
 		\Psi(t)\leq \Psi(1)=1\;\;\mbox{for}\;\; |x|\leq \dfrac{m}{2+m}.
 	\end{align*}
 	This proves the desired inequality.\vspace{2mm}
 	
 	To show the constant ${m}/{(2+m)}$ is best possible, given $a\in (0,1),$ we consider slice regular function $f_a$ given by \eqref{Eq-1.1}. For $f_a,$ a straightforward computation yields that
 	\begin{align*}
 		\mathcal{C}_{f_a}(x)&=a^m+(1-a^2)\sum_{k=1}^{\infty}|x|^ka^{k-1}+\left(\dfrac{1}{1+a}+\dfrac{|x|}{1-|x|}\right)(1-a^2)^2\sum_{k=1}^{\infty}|x|^{2k}a^{2(k-1)}\\&= a^m+\dfrac{(1-a^2)|x|}{1-a|x|}>1\;\;\mbox{for}\;\;|x|> \dfrac{m}{2+m}=R_m.
 	\end{align*}
 	This establishes that the constant \( R_m = \frac{m}{2 + m} \) is the best possible, concluding the proof.
 \end{proof}
 
 \begin{proof}[\bf Proof of Theorem \ref{BS-thm-1.3}]
  Considering \eqref{Eq-1.2} and \eqref{Eq-2.4}, we derive
 	\begin{align}\label{Eq-1.5}
 		S^*_{x}\leq (1-|a_0|^2)^2\sum_{k=1}^{\infty}k|x|^2=(1-|a_0|^2)^2\dfrac{|x|^2}{\left(1-|x|^2\right)^2}\;\;\mbox{for}\;\;|x|<1.
 	\end{align}
 	
 	Let $|a_0|=t.$  Using \eqref{Eq-1.3}, \eqref{Eq-1.5}, \eqref{Eq-2.4}, and \eqref{Eq-2.5}, a simple computation yields the following inequality
 	\begin{align}\label{Eq-33.6}
 		\mathcal{D}_f(x)&\leq 1+ m(t-1)+ \left(1-t^2\right)\frac{|x|}{1-|x|}+\sum_{k=1}^{N}d_k\left(\dfrac{(1-t^2)|x|}{1-|x|^2}\right)^{2k}\\&= 1+\Phi(t,|x|),\nonumber
 	\end{align}
 	where \begin{align*}
 		\Phi(t,|x|):=\dfrac{(1-t^2)|x|}{1-|x|}+\sum_{k=1}^{N}d_k\left(\dfrac{(1-t^2)|x|}{1-|x|^2}\right)^{2k}-m(1-t).
 	\end{align*}
 	For all \( t \in [0, 1) \), a straightforward calculation shows that \( \Phi(t, |x|) \) is a monotonically increasing function of \( |x| \). Consequently, for \( t \in [0, 1) \), we obtain $\Phi(t,|x|)\leq \Phi(t,m/(2+m)).$ By a straightforward calculation, we see that
 	\begin{align*}
 		\Phi(t,m/(2+m))=\dfrac{(1-t^2)}{2}\left(m+2F_N(t)-\dfrac{2m}{1+t}\right)=\dfrac{(1-t^2)}{2}\mathcal{W}(t),
 	\end{align*}
 	where \begin{align*}
 		F_N(t):=\sum_{k=1}^{N}d_k(1-t^2)^{2k-1}M^{2k}_m\;\;\mbox{and}\;\; \mathcal{W}(t):=m+2F_N(t)-\dfrac{2m}{1+t}.
 	\end{align*}
 	From \eqref{Eq-33.6}, we see that the desired inequality will hold if we show that \( \Phi(t, |x|) \leq 0 \). To do so, it suffices to prove that \( \mathcal{W}(t) \leq 0 \) for \( t \in [0, 1] \). A straightforward calculation shows that for \( t \in [0, 1] \), we have
 	\begin{align*}
 		t(1+t)^2M_m^2&\leq 4 M_m^2,\\ t(1+t)^2(1-t^2)^2M_m^4&\leq c_2M_m^4,\\ \vdots \\t(1+t)^2(1-t^2)^{2N-2}M_m^{2N}&\leq c_NM_m^{2N}.
 	\end{align*}
 	This leads us to the conclusion that
 	\begin{align*}
 		\mathcal{W}^{\prime}(t)&=\frac{2}{(1+t)^2}\bigg(m-2d_1t(1+t)^2M_m^{2}-6d_2t(1+t)^2(1-t^2)^2M_m^{4}-\cdots\\&\quad-2(2N-1)d_Nt(1+t)^2(1-t^2)^{2N-2}M_m^{2N}\\&\geq \dfrac{2}{(1+t)^2} (m-L(d_1,\cdots,d_N))\geq 0,
 	\end{align*}
 	if  the coefficients $d_i$ of the polynomial $Q_N$ satisfy the condition outlined in $(\ref{Eq-1.6})$. This means that $\mathcal{W}(t)$ behaves as an increasing function in $t\in [0,1]$, leading to the conclusion that $\mathcal{W}(t)\leq\mathcal{W}(1)=0$. This, in turn, establishes the desired inequality.\vspace{1.2mm}
 	
 	To show the sharpness, given $a\in (0,1),$ we consider slice regular function $f_a$ given by \eqref{Eq-1.1}. In view of \eqref{Eq-1.2} and \eqref{Eq-1.3}, we obtain
 	\begin{align*}
 		\mathcal{D}_{f_a}(x)=1-(1-a)\Psi^*(|x|,a),
 	\end{align*}
 	where 
 	\begin{align*}
 		\Psi^*(|x|,a)&:=\dfrac{1-|x|^m}{1-|x|}-\dfrac{(1+a)|x|}{1-|x|}-\dfrac{d_1|x|^2(1-a)(1+a)^2}{(1-a^2|x|^2)^2}-\dots\\&\quad-\dfrac{d_N|x|^{2N}(1-a)^{2N-1}(1+a)^{2N}}{(1-a^2|x|^2)^{2N}}.
 	\end{align*}
 	Consequently, we arrive at the conclusion that
 	\begin{align}
 		\lim\limits_{a\to 1^{-}}\Psi^*(|x|,a)=m-\frac{2|x|}{1-|x|}<0\;\mbox{for}\;\;|x|>R_m=\dfrac{m}{2+m}.
 	\end{align}
 	Thus, it follows that $\Psi^*(|x|,a)<0$ for $a$ sufficiently close to $1$. 
 	Hence, we have
 	\begin{align*}
 		\mathcal{D}_{f_a}(x)=1-(1-a)\Psi^*(|x|)>1\; \mbox{for}\;|x|>R_m
 	\end{align*}
 	which shows that the constant $R_m={m}/{(2+m)}$ is best possible and thereby concluding the proof.
 \end{proof}

 \begin{proof}[\bf Proof of Theorem \ref{Thm-1.7}]
 	Let $\alpha=1-a_0$. Under the condition of Theorem \ref{Thm-1.7}, it follows that  (see\cite[p.1231]{Ren-Wang-CAOT-2015})
 	\begin{align}\label{Lem-2.11}
 		|a_k|\leq 2\left(1-{\rm Re}(f(0))\right)\;\;\;\mbox{for}\;\;\; k\in\mathbb{N}.
 	\end{align}
 	In view of \eqref{Lem-2.11}, we obtain $|a_k|\leq 2\alpha$ for all $k\geq 1$. For $\alpha\in (0, 1]$ and $|x|=r<1,$ we obtain 
 	\begin{align}\label{Eq-33.9}
 		\mathcal{E}_f(x)&\leq 1-\alpha+\frac{2\alpha r}{1-r}+\left(\frac{1}{1+a_0}+\frac{r}{1-r}\right)\frac{4\alpha^2r^2}{1-r^2}\\&\nonumber=1-\alpha\left(\frac{1-3r}{1-r}-\left(\frac{1+r-\alpha r}{(1-r)(2-\alpha)}\right)\frac{4\alpha r^2}{1-r^2}\right)\\&\nonumber=1-\alpha\frac{Q(\alpha, r)}{(2-\alpha)(1-r)(1-r^2)},
 	\end{align}
 	where \[ Q(\alpha, r):=4r^3\alpha^2-(7r^3+3r^2-3r+1)\alpha+6r^3-2r^2-6r+2.\]
 Our aim is to show that $ Q(\alpha, r)\geq 0$ for $r\leq R_*$ and $\alpha\in (0, 1]$. It follows that $\frac{\partial^2 Q(\alpha, r)}{\partial \alpha^2}\geq 0$ for any $\alpha\in (0, 1]$ and thus, $\frac{\partial Q}{\partial \alpha}$ is an increasing function of $\alpha$. This gives
 	\begin{align*}
 		\frac{\partial Q(\alpha, r)}{\partial \alpha}\leq \frac{\partial Q(1, r)}{\partial\alpha}=-(1-r)^3<0\;\;\mbox{for}\;\;r<1.
 	\end{align*}
This shows that $Q$ is a decreasing function of $\alpha$ on $(0, 1]$ so that $Q(\alpha, r)\geq Q(1, r)=3r^3-5r^2-3r+1$	which is greater than or equal to $0$ for all $r\leq R_*$, where $R_*\approx 0.24683$ is the unique root of the equation $3r^3-5r^2-3r+1=0$ which lies in $(0, 1)$. It is clear from \eqref{Eq-33.9} that $\mathcal{E}_f(x)\leq 1$ for $r\leq R_*$ which is the desired inequality. \vspace{1.2mm}
 	
 	To prove the sharpness, for $u\in \partial \mathbb{B}$ and $a\in (0,1),$ we consider the function $g_a$ which is defined by
 
 \begin{align}\label{Eq-3.6}
 	g_a(x):=a-2(1-a)x(1-x)^{-\bullet}u=a-\sum_{k=1}^{\infty}x^ku2(1-a).
 \end{align}
 
 	For $g_a,$ we obtain 
 	\begin{align*}
 		\mathcal{E}_{g_a}(x)&=\sum_{k=0}^{\infty}|x^ka_k|+\left(\frac{1}{1+a_0}+\frac{|x|}{1-|x|}\right)\sum_{k=1}^{\infty}|x^{k}a_k|^2\\&=a+2(1-a)\dfrac{|x|}{1-|x|}+\left(\frac{1}{1+a}+\frac{|x|}{1-|x|}\right)\dfrac{4(1-a)^2|x|^2}{1-|x|^2}.
 	\end{align*}
 By the similar argument used above, it can be easily shown that $\mathcal{E}_{g_a}(x)>1$ for $|x|>R_*.$ This completes the proof.
 \end{proof}
 
 \begin{proof}[\bf Proof of Theorem \ref{Theom-1.7}]
 Clearly, $\mathcal{F}_{f}$ is an increasing function of $|x|=r$ and hence, we only need to prove the inequality $\mathcal{F}_{f}(|x|)\leq 1$ for $|x|=\xi={1}/{(5-2a_0)}$. In view of \eqref{Lem-2.11}, we obtain the inequality
 
\begin{align*}
	 S^*_{x}=\sum_{k=1}^\infty k|x^ka_k|^2 \leq 4(1-a_0)^2\frac{r^2}{(1-r)^2},\;\;|x|=r.
\end{align*}
 
In light of \eqref{Lem-2.11}, we obtain
 \begin{align*}
 	\mathcal{F}_f(\xi)&= a_0+\frac{1-a_0}{2-a_0}+\left(\frac{1}{1+a_0}+\frac{1}{4-2a_0}\right)\frac{(1-a_0)^2}{(2-a_0)(3-a_0)}\\&\quad+\beta\frac{(1-a_0)^2(5-2a_0)^2}{4(2-a_0)^2(3-a_0)^2}\\&=1+\frac{(1-a_0)}{4(2-a_0)^2(3-a_0)^2(1+a_0)}G(a_0),
 \end{align*}
 where 
 \begin{align}\label{Eq-1.100}
 	G(t):=-4t^5+32t^4-82t^3+58t^2+38t-42-4\beta t^4+20\beta t^3-21\beta t^2-20\beta t+25\beta.
 \end{align}
 To prove the desired inequality, it is sufficient to show that $G(t)\leq 0$ for all $t\in [0, 1]$. We see that $G(1)=0$ and $G(0)=-42+25\beta\leq 0$ for $\beta\leq 1.68$. Next, we find that 
 \begin{align*}
 	G^{\prime}(t)=-20t^4+128t^3-246 t^2+116t+38-16\beta t^3+60\beta t^2-42\beta t-20\beta
 \end{align*}
 and therefore, $ G^{\prime}(1) =16-18\beta$ showing that $G^{\prime}(1)\geq 0$ for $\beta\leq \frac{8}{9}$ and $G^{\prime}(1)\leq 0$ for $\beta\geq \frac{8}{9}$. Since $G(1)=0$ and $G$ is a differentiable function, so that if $G^{\prime}(1)<0$, then there exists a $t_0\in (1-\epsilon, 1)$ such that $G(t_0)>0$ for some $\epsilon>0$, which does not meet the requirements.  Therefore, we should satisfy $G^{\prime}(1)\geq 0$, that is $\beta\leq \frac{8}{9}$. When $\beta\in [0, 8/9]$, we have $G(0)=-42+25\beta\leq 0$ and $G(1)=0$. Since
 \begin{align*}
 	\frac{\partial G^{\prime}}{\partial \beta}=-16t^3+60 t^2-42t-20\leq 0\; \mbox{for}\; t\in [0, 1],
 \end{align*}
 so that for $\beta\in [0, 8/9]$, we have 
 \begin{align*}
 	G^{\prime}(t)&\geq -20 t^4+128t^3-\frac{128}{9}t^3-246 t^2+\frac{160}{3}t^2+116t -\frac{112}{3}t+38-\frac{160}{9}:=L(t).
 \end{align*}
 We see that $L(0)\approx 20.2222$ and $L(1)=0$, and with the aid of the Mathematica, we see that $L(t)=0$ has no root in $[0, 1)$ and $L(t)$ is continuous. Then $G^{\prime}(t)\geq L(t)>0$ for all $t\in [0, 1)$. This implies that $G(t)\leq G(1)=0$ in $[0, 1]$. That is $\mathcal{F}_{f}\left(1/(5-2a_0)\right)\leq 1$ for $\beta\in [0, 8/9]$. In conclusion, the maximum value that $\beta$ can take as $8/9$. Thus $G(t)\leq 0$ for $\beta=8/9$ and hence for all $0<\beta\leq 8/9$.\vspace{1.2mm}
 
 To prove  that the constant $\beta=8/9$ is sharp, for some $u\in \partial \mathbb{B}$ and $a\in (0,1),$ we consider the function $g_a$ given in \eqref{Eq-3.6}. For the function $g_a$, we compute for $x=1/(5-2a),$
 \begin{align*}
 	\mathcal{F}_{g_a}(x)=\sum_{k=0}^{\infty}|x^ka_k|+\left(\frac{1}{1+a_0}+\frac{|x|}{1-|x|}\right)\sum_{k=1}^{\infty}|x^ka_k|^2+\beta_1 S^*_{x},
 \end{align*}
 which simplifies to
 \[ \mathcal{F}_{g_a}(x)=1+\dfrac{(1-a)	G(a)}{4(1-a)^2(3-a)^2(1+a)},\]
 where $G(a)$ is given in \eqref{Eq-1.100}. For $\beta_1>\beta=8/9,$ we get that $G^\prime(1)<0$ and $G(0.999999)>0,$ which implies $\mathcal{F}_f(x)>1.$ This proves the sharpness assertion, and with this, the proof is complete.
 \end{proof}\vspace{1.2mm}
 		
 \noindent{\bf Acknowledgment:} The research of the first author is supported by UGC-JRF (Ref. No. 201610135853), New Delhi, Govt. of India and second author is supported by SERB File No. SUR/2022/002244, Govt. of India.\vspace{2mm}

 
 \noindent\textbf{Conflict of interest:} The authors declare that there is no conflict  of interest regarding the publication of this paper.\vspace{1.5mm}
 
 \noindent\textbf{Data availability statement:}  Data sharing not applicable to this article as no datasets were generated or analysed during the current study.\vspace{1.5mm}
 
 
 \noindent {\bf Authors' contributions:} Both the authors have equal contributions.
  
\end{document}